\documentclass[a4paper,10pt]{amsart}
\usepackage{amsmath,amsfonts,amssymb, amsthm, xypic}
\usepackage{epsfig,psfrag}
\usepackage[all]{xy}

\def\ss{\textbf{ss}}

\def\kF{\mathcal{F}}
\def\kG{\mathcal{G}}
\def\kL{\mathcal{L}}

\def\Z{\mathbb{Z}}

\def\E{{X}}

\def\H{\mathbf{H}}

\def\CC{\mathbb{C}}

\def\U{\boldsymbol{\mathcal{E}}}

\def\qed{\hfill $\sqcap \hskip-6.5pt \sqcup$}

\def\a{\alpha}

\def\qlb{\overline{\mathbb{Q}}_l}

\def\C{\mathbb{C}}
\def\UU{\mathbf{U}}
\def\U{\mathbf{U}}

\def\fqb{\overline{\mathbb{F}_q}}
\def\fq{\mathbb{F}_q}

\newtheorem{theo}{\bf{Theorem}}[section]
\newtheorem{lem}[theo]{Lemma}

\newtheorem{cor}[theo]{Corollary}

\newtheorem{prop}[theo]{Proposition}

\numberwithin{equation}{section}

\title{Spherical Hall algebras of curves and Harder-Narasimhan stratas}

\author{O. Schiffmann}

\begin{document}
\maketitle

\begin{abstract}
We show that the characteristic function $\mathbf{1}_{S_{\underline{\a}}}$ of any Harder-Narasimhan strata $S_{\underline{\a}} \subset Coh^{\a}_X$ belongs to the spherical Hall algebra $\H^{sph}_X$ of a smooth projective curve $X$ (defined over a finite field $\mathbb{F}_q$). We prove a similar result in the geometric setting : the intersection cohomology complex $IC(\underline{S}_{\underline{\a}})$ of any Harder-Narasimhan strata $\underline{S}_{\underline{\a}} \subset \underline{Coh}^{\a}_{\overline{X}}$ belongs to the category $\mathcal{Q}_X$ of spherical Eisenstein sheaves of $X$. We show by a simple example how a complete description of all spherical Eisenstein sheaves would necessary involve the Brill-Noether stratas of $\underline{Coh}^{\a}_{\overline{X}}$.
\end{abstract}

\vspace{.2in}

{{\centerline{JOYEUX ANNIVERSAIRE GERARD !}}}

\vspace{.2in}

\setcounter{section}{-1}

\section{Introduction}

\vspace{.1in}

Let $X$ be a smooth projective curve defined over a finite field $\fq$. To such a curve is associated, through a general formalism developed by Ringel, a Hopf algebra $\H_X$ (called the \textit{Hall algebra of $X$}). As a vector space, $\H_X$ consists of all finitely supported functions on the set of (isomorphism classes of) coherent sheaves over $X$, and the (co)product encodes the structure of the extensions between coherent sheaves over $X$ (see e.g. \cite{SLectures}).

\vspace{.1in}

Hall algebras were first considered by Ringel in the context of representations of quivers. He showed that a certain natural subalgebra $\mathbf{C}_{\vec{Q}}$ of the Hall algebra $\H_{\vec{Q}}$ of a quiver $\vec{Q}$ is isomorphic
to the quantized enveloping algebra $\UU^+_q(\mathfrak{g})$ of a Kac-Moody algebra $\mathfrak{g}$ attached to $\vec{Q}$ (\cite{Ringel}). This discovery paved the way for a completely new approach to the theory of quantum groups based on the representation theory of quivers (see e.g.  \cite{KS}, \cite{Nak}, \cite{VV},...)
One of the most important development is the work of Lusztig who considered a geometric version of the Hall algebra of a quiver, in which functions are replaced by
constructible sheaves (on moduli stacks); this gave rise to the theory \textit{canonical bases} of quantum groups, whose impact in algebraic and geometric representation theory is well-known (see \cite{Lusztig}).

\vspace{.1in}

In the context of smooth projective curves, Hall algebras first appeared in pioneering work of Kapranov, in relation to the theory of automorphic forms over function fields (see \cite{Kap}). He observed some striking similarities between the Hall algebras of curves and quantum affine algebras $\U^+_q(\widehat{\mathfrak{g}})$ (more precisely between the \textit{functional equations for Eisenstein series} over function fields and the so-called \textit{Drinfeld relations} in quantum affine algebras). This analogy becomes very precise for $X=\mathbb{P}^1$ (see \cite{BK}). Motivated by the theory of quantum loop algebras, we introduced in \cite{SDuke} a natural subalgebra $\H_X^{sph}$ of $\H_X$ which we call the \textit{spherical Hall algebra of $X$}; in the language of automorphic forms, $\H_X^{sph}$ is generated by the Fourier coefficients of all Eisenstein series induced from the trivial character of a maximal torus. Inspired by work of Laumon, we also singled out a category $\mathcal{Q}_X$ of simple perverse sheaves over the moduli stacks $\underline{Coh}^{r,d}_X$ of coherent sheaves over $\overline{X}=X \otimes \fqb$. These simple perverse sheaves (which we call \textit{spherical Eisenstein sheaves}) are expected to provide (by means of the Faisceaux-Fonctions correspondence) a canonical basis for the spherical Hall algebras $\H^{sph}_X$. From \cite{SV3} it is also natural to expect $\mathcal{Q}_X$ to play a role in the geometric Langlands program for local systems in the formal neighborhood of the trivial local system.

\vspace{.1in}

The spherical Hall algebras $\H_X$ and the spherical Eisenstein sheaves $\mathcal{Q}_X$ of several low-genus (possibly orbifold) curves were computed in a series of papers (see \cite{SDuke}, \cite{BS}, \cite{SInvent}, \cite{Scano}) where they were shown to yield interesting quantum loop algebras (such as quantum toroidal algebras or spherical Cherednik algebras of type $A$) equipped with some canonical bases. Much less is known for higher genus curves; a combinatorial realization of $\H_X^{sph}$ as a shuffle algebra is given in \cite{SV3} for an arbitrary curve, but it is rather hard to analyze directly algebraically.

\vspace{.1in}

In this note, as a first step towards understanding these higher genus spherical Hall algebras, we exhibit an explicit class of elements in $\H_X^{sph}$ when $X$ is of genus $g >1$.  Namely we prove (see Theorem~\ref{T:sstable}) that the characteristic functions of all the Harder-Narasimhan stratas $S_{\underline{\a}}$ belong to $\H_X^{sph}$ (see Section~2 for notations). We also give a geometric version of the same result : the intersection cohomology complex $IC(\underline{S}_{\underline{\a}})$ of any Harder-Narasimhan strata is a spherical Eisenstein sheaf (Theorem~\ref{T:maingeom}). As we show by an example, these classes of functions (resp. simple perverse sheaves) come very far from exhausting the whole of $\H_X^{sph}$ (resp. $\mathcal{Q}_X$) : a full description would (at least) involve the various Brill-Noether loci in $\underline{Coh}^{r,d}_X$ (see Remark~5.3.).

\vspace{.2in}

\section{Spherical Hall algebras of curves}

\vspace{.15in}

\paragraph{\textbf{1.1}} Let $X$ be a connected smooth projective curve of genus $g$ defined over the finite field $\fq$. We will assume here that $g >1$, although most of the results proved in this note hold for rational and elliptic curves as well (see \cite{SInvent}, \cite{Scano}). Let $Coh(X)$ stand for the category of coherent sheaves over $X$. Let us denote by $\langle \;,\;\rangle$ the Euler form on the Grothendieck group $K_0(X)$, and let $K'_0(X)=K_0(X)/ \text{rad} \langle\;,\;\rangle$ be the numerical Grothendieck group of $X$. We have $K'_0(X)=\Z \;\text{rk} \oplus \Z \;\text{deg}$ where $\text{rk}$ and $\text{deg}$ are the rank and degree functions respectively. The Euler form on $K_0(X)$ is given by
\begin{equation}\label{E:Euler1}
\langle \kF, \kG \rangle = (1-g) \text{rk}(\kF) \text{rk}(\kG) + \left| \begin{array}{cc} \text{rk}(\mathcal{F}) & \text{rk}(\mathcal{G}) \\
\text{deg}(\mathcal{F}) & \text{deg}(\mathcal{G}) \end{array} \right|.
\end{equation}

\vspace{.2in}

\paragraph{\textbf{1.2.}} Let us briefly recall the definition of the Hall algebra $\H_X$ of $X$. We refer the reader to \cite[Lecture~4]{SLectures} for more details. Let $\mathcal{I}_X$ stand for the set of all coherent sheaves over $X$. Put
$$\H_X=\{f: \mathcal{I}_X \to \C\;|\; \# \;supp\;f < \infty\}=\bigoplus_{ \mathcal{F} \in \mathcal{I}} \C 1_{\mathcal{F}}$$
where $1_{\mathcal{F}}$ is the characteristic function of the point $\mathcal{F} \in \mathcal{I}$.
Let us fix a square root $v$ of $q^{-1}$.
The multiplication in $\H_X$ is defined by the following formula
$$(f \cdot g) (\mathcal{R})=\sum_{\mathcal{N} \subseteq
\mathcal{R}} v^{-\langle \mathcal{R}/\mathcal{N},\mathcal{N}\rangle} f(\mathcal{R}/\mathcal{N}) g(\mathcal{N})$$
and the comultiplication is
$$\Delta(f) (\mathcal{M},\mathcal{N})=\frac{v^{\langle \mathcal{M},\mathcal{N}\rangle}}
{|\text{Ext}^1(\mathcal{M},\mathcal{N})|}\sum_{\xi \in
\text{Ext}^1(\mathcal{M},\mathcal{N})} f(\mathcal{X}_{\xi})$$ where
$\mathcal{X}_{\xi}$ is the extension of $\mathcal{N}$ by
$\mathcal{M}$ corresponding to $\xi$. Notice that the coproduct
$\Delta$ takes values in a completion $\H_{X}
\widehat{\otimes} \H_{X}$ of the tensor space  $\H_{X}
\otimes \H_{X}$ only (see e.g. \cite[Section~2]{BS}). The triple $(\H_{\E},\cdot, \Delta)$ is not a
(topological) bialgebra, but it becomes one if we suitably twist the
coproduct. For this we introduce an extra subalgebra
$\boldsymbol{\mathcal{K}}=\CC[\boldsymbol{\kappa}_{r,d}]$, $(r,d)
\in \Z^2$, and we define an extended Hall algebra
$\widetilde{\H}_{\E} = \H_{\E} \otimes \boldsymbol{\mathcal{K}}$
with relations
$$\boldsymbol{\kappa}_{r,d}\, \boldsymbol{\kappa}_{s,l}=\boldsymbol{\kappa}_{r+s,d+l},
\quad \boldsymbol{\kappa}_{0,0}=1,\quad \boldsymbol{\kappa}_{r,d}
1_{\mathcal{M}}\, \boldsymbol{\kappa}_{r,d}^{-1} =
v^{-2r(1-g)\text{rk}(\mathcal{M})} 1_{\mathcal{M}}.$$ 
The new coproduct is given by the formulas
$$\widetilde{\Delta}(\boldsymbol{\kappa}_{r,d})=\boldsymbol{\kappa}_{r,d} \otimes \boldsymbol{\kappa}_{r,d},$$
$$\widetilde{\Delta}(f)=\sum_{\mathcal{M},\mathcal{N}} \Delta(f)(\mathcal{M},\mathcal{N})
1_{\mathcal{M}}\boldsymbol{\kappa}_{r_{\mathcal{N}},d_{\mathcal{N}}}
\otimes 1_{\mathcal{N}}.$$ Then $(\widetilde{\H}_\E, \cdot,
\widetilde{\Delta})$ is a topological bialgebra. 

\vspace{.1in}

The Hall algebras $\H_{\E}$ and $\widetilde{\H}_{\E}$ are $\Z^2$-graded (by the class in the numerical Grothendieck group).
We will sometimes write $\Delta_{\a,\beta}$ or $\widetilde{\Delta}_{\a,\beta}$ in order to specify the
graded components of the coproduct.

\vspace{.2in}

\paragraph{\textbf{1.3.}} We will especially be interested in the \textit{spherical} subalgebra $\H^{sph}_X$ of $\H_X$, which is defined as follows. For any $d \in \Z$ let $Pic^d(X)$ stand for the (finite) set of line bundles over $X$ of degree $d$, and let us set
$$\mathbf{1}^{\ss}_{1,d} =\sum_{\mathcal{L} \in Pic^d(X)} 1_{\mathcal{L}}.$$
Next, let $d \geq 1$ and let $Tor^d(X)$ stand for the (finite) set of all torsion sheaves over $X$ of degree $d$. We now set
$$\mathbf{1}_{0,d}=\sum_{\mathcal{T} \in Tor^d(X)} 1_{\mathcal{T}}.$$
The spherical Hall algebra $\H_X^{sph}$ is generated by the elements $\{\mathbf{1}^{\ss}_{1,d}\;|\; d \in \Z\} \cup \{ \mathbf{1}_{0,d}\;|\; d \geq 1\}$. In an effort to unburden the notation, and because this should not cause any confusion here, we will simply write $\UU_X$ for $\H_X^{sph}$.
The spherical Hall algebra contains two natural subalgebras, namely $\UU^>_{X}$ which is generated by $\{\mathbf{1}^{\ss}_{1,d}\;|\; d \in \Z\}$, and $\UU^0_X$ which is generated by $\{ \mathbf{1}_{0,d}\;|\; d \geq 1\}$. Moreover, the multiplication map gives an isomorphism
$\UU^>_X \otimes \UU^0_X \to \UU_X$ (see e.g. \cite[Section~6]{SV}).

\vspace{.2in}

\section{Harder-Narasimhan stratifications}

\vspace{.15in}

\paragraph{\textbf{2.1.}} Let us now briefly recall the various notions related to semistablity of coherent sheaves over curves. We refer to \cite{HN}, \cite{Sh} for more details. We fix a smooth projective curve $\E$ of genus $g$ as in Section~1.1. The slope of a coherent sheaf $\mathcal{F}$ over $\E$ is defined to be
$$\mu(\mathcal{F})=\frac{deg(\mathcal{F})}{rank(\mathcal{F})} \in \mathbb{Q} \cup \{\infty\}.$$
 The sheaf $\mathcal{F}$ is said to be \textit{semistable of slope $\nu$} if $\mu(\mathcal{F})=\nu$ and if $\mu(\mathcal{G}) \leq \nu$ for any subsheaf $\mathcal{G}$ of $\mathcal{F}$. If the above condition holds with $<$ instead of $\leq$ then we say that $\mathcal{F}$ is \textit{stable}. We denote by $\textbf{C}_{\nu}$ the full subcategory of $Coh(\E)$ whose objects are semistable sheaves of slope $\nu$. The categories $\textbf{C}_\nu$ are abelian and artinian. The simple objects of $\textbf{C}_{\nu}$ are precisely given by the stable sheaves of slope $\nu$.

\vspace{.1in}

The fundamental properties of the categories $\textbf{C}_{\nu}$ are listed below.

\vspace{.1in}

\begin{prop}\label{P:A1} The following hold~:
\begin{enumerate}
\item[i)] $\text{Hom}(\textbf{C}_{\nu}, \textbf{C}_{\eta})=0$ if $\nu > \eta$,
\item[ii)] $\text{Ext}(\textbf{C}_{\nu}, \textbf{C}_{\eta})=0$ if $\eta > \nu + 2(g-1)$,
\item[iii)] any coherent sheaf $\mathcal{F}$ possesses a unique filtration
\begin{equation}\label{E:Afiltr}
0 \subsetneq \mathcal{F}_l \subsetneq \cdots \subsetneq \mathcal{F}_1=\mathcal{F}
\end{equation}
satisfying the following conditions~: $\mathcal{F}_i/\mathcal{F}_{i+1}$ is semistable for all $i$ and 
$$\mu(\mathcal{F}_1/\mathcal{F}_2) < \cdots < \mu(\mathcal{F}_{l-1}/\mathcal{F}_l) < \mu(\mathcal{F}_l).$$
\end{enumerate}
\end{prop}

\vspace{.1in}

The filtration (\ref{E:Afiltr}) is called the \textit{Harder-Narasimhan (or HN) filtration} of $\mathcal{F})$. We also define the
\textit{HN-type} of $\mathcal{F}$ to be $HN(\mathcal{F})=(\a_1, \ldots, \a_l)$ with $\a_i= \overline{\mathcal{F}_{i}}-
\overline{\mathcal{F}_{i+1}}$. Here $\overline{\mathcal{G}}=(rank(\mathcal{G}), deg(\mathcal{G})) \in \Z^2$ is the
class of a sheaf $\mathcal{G}$ in the (numerical) Grothendieck group of $Coh(\E)$-- see Section~1.1. Note that the weight $\a:=\a_1+ \cdots + \a_l$ of the HN type of $\mathcal{F}$ is equal to $\overline{\mathcal{F}}$.

\vspace{.1in}

It is convenient to view an HN type $(\a_1, \ldots, \a_l)$ as a polygon as follows~:

\vspace{.2in}

\centerline{
\begin{picture}(200,120)
\put(0,0){\line(0,1){120}}
\put(0,60){\line(1,0){200}}
\put(0,60){\circle*{2}}
\put(-7,60){$\textbf{o}$}
\put(40,30){\circle*{2}}
\put(35,24){$\alpha_1$}
\put(0,60){\line(4,-3){40}}
\put(40,30){\line(5,-2){50}}
\put(90,10){\circle*{2}}
\put(78,4){$\alpha_1+\alpha_2$}
\put(90,10){\line(2,1){46}}
\put(136,33){\circle*{2}}
\put(136,33){\line(1,2){20}}
\put(156,73){\circle*{2}}
\put(156,73){\line(1,4){8}}
\put(164,105){\circle*{2}}
\put(166,102){$\a$}
\end{picture}}
\vspace{.05in}
\centerline{\textbf{Figure 1.} A Harder-Narasimhan polygon of weight $\a$.}

\vspace{.15in}

This polygon, called the \textit{HN polygon} of $\mathcal{F}$, is convex by construction. The following useful
result is a consequence of Proposition~\ref{P:A1} (see e.g. \cite{Sh}, Theorem~2).

\vspace{.1in}

\begin{prop}\label{P:A2} Let $\mathcal{F}$ be a coherent sheaf over $\E$ of class $\a \in \Z^2$. Let $0 \subsetneq \mathcal{F}_l \subsetneq \cdots \subsetneq \mathcal{F}_1=\mathcal{F}$ be the HN filtration of $\mathcal{F}$. Let $\mathcal{G}$ be a subsheaf of $\mathcal{F}$ of class $\gamma$. Then
\begin{enumerate}
\item[i)] the point $\beta:=\alpha-\gamma$ lies above the HN polygon of $\mathcal{F}$,
\item[ii)] if moreover $\beta$ is a vertex of the HN polygon of $\mathcal{F}$, that is if $\gamma=\a_i+ \cdots + \alpha_l$
for some $1 \leq i \leq l$, then $\mathcal{G}=\mathcal{F}_i$.
\end{enumerate}
\end{prop}

\vspace{.2in}

\paragraph{\textbf{2.2.}} We may stratify the set $\mathcal{I}_X$ of all isomorphism classes of coherent sheaves over $\E$ by the HN-type and write $\mathcal{I}_X=\bigsqcup_{\underline{\a}} S_{\underline{\a}}$ where $\underline{\a}$ runs through the set of all possible
HN types, i.e. tuples $\underline{\a}=(\a_1, \ldots, \a_l)$ with $\a_i \in (\Z^2)^+$ and $\mu(\a_1) < \cdots < \mu(\a_l)$. Here $(\Z^2)^+=\{ (r,d) \in \Z^2\;|\; r \geq 1\;or\;r=0, d >0\}$. For instance, if $\underline{\a}=(\a)$ then $S_{\underline{\a}}$ is the set of isomorphism classes of semistable sheaves of class $\a$.  Let us denote by $\mathbf{1}_{S_{\underline{\a}}} \in \H_{\E}$ the characteristic function of the set of sheaves of a fixed HN type $\underline{\a}$. Since $\E$ is defined over a finite field, $S_{\underline{\a}}$ is finite for any $\underline{\a}$ hence $\mathbf{1}_{S_{\underline{\a}}}$ is a well-defined element of $\H_{\E}$. For $\a \in (\Z^2)^+$ we will simply denote by $\mathbf{1}^{\ss}_{\a}$ the characteristic function of $S_{\alpha}$.

\vspace{.1in}

From the uniqueness of the HN filtration of a given coherent sheaf we easily deduce

\vspace{.1in}

\begin{prop}\label{P:A3} For any HN type $\underline{\a}=(\a_1, \ldots, \a_l)$ we have
$$\mathbf{1}_{S_{\underline{\a}}}=v^{\sum_{i<j} \langle \a_i, \a_j\rangle} \mathbf{1}^{\ss}_{\a_1} \cdots \mathbf{1}^{\ss}_{\a_l}.$$
\end{prop} 

\vspace{.1in}

We use the stratification by HN type to define a completion of the Hall algebra $\H_{\E}$ as follows. For $n \in \Z$ let us write $\underline{\a} \geq n$ if $\underline{\a}=(\a_1, \ldots, \a_l)$ with $\mu(\a_1) \geq n$. Let $\mathbf{C}_{\geq n}$ be the full subcategory of $Coh(\E)$ generated by $\mathbf{C}_{\nu}$ for all $\nu \geq n$. By definition, the HN type $\underline{\a}$ of a sheaf $\mathcal{F}$ satisfies $\underline{\a} \geq n$ if and only if $\mathcal{F} \in \mathbf{C}_{\geq n}$.

 The set of HN types of a fixed weight $\a$ satisfying $\underline{\a} \geq n$ is finite for any $n$. Let $\H_{\E}^{< n}[\a]$ be the subspace of $\H_{\E}[\a]$ consisting of functions supported on the \textit{complement} of $\bigcup_{\underline{\a} \geq n} S_{\underline{\a}}$. It is a subspace of $\H_{\E}[\a]$ of finite codimension. Moreover there are some obvious inclusions 
$\H^{< m}_{\E}[\a] \to \H^{< n}_{\E}[\a]$ for any $m<n$. Put $\H^{\geq n}_{\E}[\a]=\H_{\E}[\a]/ \H^{<n}_{\E}[\a]$.  This is a finite dimensional space. We put
$$\widehat{\H}_{\E}[\a]=\underset{\longrightarrow}{\text{Lim}}\; \H^{\geq n}_{\E}[\a], \qquad \widehat{\H}_{\E}=\bigoplus_{\a} \widehat{\H}_{\E}[\a].$$
Note that $\widehat{\H}_{\E}[\a]=\{f~: \mathcal{I}_{\a} \to \CC\}=\prod_{\mathcal{F} \in \mathcal{I}_{\a}} \CC 1_{\mathcal{F}}$ as a vector space, where we have denoted by $\mathcal{I}_{\a} \subset \mathcal{I}_X$ the set of all coherent sheaves of class $\a$. It is shown in \cite{BS}, Section~2 that the product and coproduct are well-defined in the limit and endow $\widehat{\H}_{\E}$ with the structure of a (twisted) bialgebra.

\vspace{.1in}

Consider the elements
$$\mathbf{1}_{\a}=\sum_{\mathcal{F} \in \mathcal{I}_{\a}} 1_{\mathcal{F}}, \qquad \mathbf{1}^{\textbf{vec}}_{\a}=\sum_{\mathcal{V} \in \mathcal{I}_{\a}^{vec}} 1_{\mathcal{V}}$$
where the second sum ranges over all (isomorphism classes of) vector bundles of class $\a$. These are both elements of $\widehat{\H}_{\E}$. As a direct corollary
of Proposition~\ref{P:A3} we have the following identities~:
\begin{equation}\label{E:A21}
\mathbf{1}_{\a}=\sum_{\underline{\a} \in X_{\a}} v^{\sum_{i < j} \langle \a_i, \a_j \rangle} \mathbf{1}^{\ss}_{\a_1} \cdots \mathbf{1}^{\ss}_{\a_l}, \qquad \mathbf{1}^{\textbf{vec}}_{\a}=\sum_{\underline{\a} \in Y_{\a}} v^{\sum_{i < j} \langle \a_i, \a_j \rangle} \mathbf{1}^{\ss}_{\a_1} \cdots \mathbf{1}^{\ss}_{\a_l},
\end{equation} 
where $X_{\a}$ is the set of all HN types of weight $\a$ and $Y_{\a}$ is the set of all HN types $\underline{\a}=(\a_1, \ldots, \a_l)$ of weight $\a$ for which $\mu(\a_l) < \infty$.

\vspace{.2in}

\section{Characteristic functions of semistables}

\vspace{.2in}

\paragraph{\textbf{3.1.}} Our aim in this section is to prove the following theorem~:

\vspace{.1in}

\begin{theo}\label{T:sstable}
For any $\a \in (\Z^2)^+$ we have $\mathbf{1}^{\ss}_{\a} \in \UU_{\E}$.
\end{theo}

\vspace{.1in}

Our proof of Theorem~\ref{T:sstable} hinges on the following preliminary result. Let us denote by $\widehat{\UU}_{\E}$ the
completion of $\UU_{\E}$ in $\widehat{\H}_{\E}$ (i.e. $\widehat{\UU}_{\E}[\a]=\underset{\longrightarrow}{\text{Lim}}\; \UU_{\E}[\a]/ ( \UU_{\E}[\a] \cap \H^{< n}_{\E}[\a])$). 

\vspace{.1in}

\begin{prop}\label{P:A4} For any $\a \in (\Z^2)^+$ we have $\mathbf{1}^{\ss}_{\a} \in \widehat{\UU}_{\E}$.
\end{prop}
\begin{proof} We may use Reineke's inversion formula (see \cite{Reineke}, Section~5.) to write
$$\mathbf{1}^{\ss}_{\a}=\sum_{\underline{\beta}} (-1)^{s-1} v^{\sum_{i<j}\langle \beta_i, \beta_j\rangle} \mathbf{1}_{\beta_1}
\cdots \mathbf{1}_{\beta_s}$$
where the sum ranges over all tuples $\underline{\beta}=(\beta_1, \ldots, \beta_s)$ of elements of $(\Z^2)^+$ satisfying
$\mu(\sum_{l=k}^s \beta_l) > \mu(\alpha)$ for all $k=1, \ldots, s$. The above sum converges in $\widehat{\H}_{\E}$. Since $\widehat{\UU}_{\E}$ is a subalgebra of $\widehat{\H}_{\E}$, the proposition will be proved if we can show that $\mathbf{1}_{\a} \in \widehat{\UU}_{\E}$ for all $\a$. Furthermore, because $\mathbf{1}_{\a}=\sum_{l \geq 0} v^{l rank(\a)} \mathbf{1}^{\textbf{vec}}_{\a-(0,l)} \mathbf{1}_{(0,l)}$
and $\mathbf{1}_{(0,l)} \in \UU_{\E}$ for all $l$, it suffices in fact to prove that $\mathbf{1}^{\textbf{vec}}_{\a} \in \widehat{\UU}_{\E}$ for all $\a$.

Let us write $\a=(r,d)$ and argue by induction on the rank $r$. The cases $r=0,1$ are obvious so
let $r>1$ and let us assume that 
\begin{equation}\label{E:Aproof1}
\mathbf{1}_{(r',d)} \in \widehat{\UU}_{\E}, \qquad \mathbf{1}^{\textbf{vec}}_{(r',d)} \in \widehat{\UU}_{\E}
\end{equation}
 for all $r'<r$. We have to show that for any $d \in \Z$ and any $n \in \Z$ it holds
\begin{equation}\label{E:huyt}
\mathbf{1}^{\textbf{vec}}_{r,d} \in \UU_{\E} + \H_{\E}^{<n}.
\end{equation}
Let us fix $n$ and argue by induction on $d$. If $d < nr$ then no vector bundle of rank $r$ and degree $d$ may belong to $\mathbf{C}_{\geq n}$ and hence have an HN type $\underline{\a} \geq n$. Therefore $\mathbf{1}^{\textbf{vec}}_{r,d} \in \H_{\E}^{< n}$. Now let us fix some $d$ and assume that (\ref{E:huyt}) holds for all $d'<d$. 

Choose $N < n-2(g-1)$ and let us consider the product $\mathbf{1}_{r-1,d-N}\cdot\mathbf{1}^{\textbf{vec}}_{1,N}$. By definition, we have
$$\mathbf{1}_{r-1,d-N}\cdot\mathbf{1}^{\textbf{vec}}_{1,N}=\sum_{\mathcal{F}} c_{\kF}[\kF]$$
where 
\begin{equation*}
c_{\kF}=v^{-\langle (r-1,d-N),(1,N)\rangle}\sum_{\kL \in Pic^N(X)}\frac{\#\{\kL \hookrightarrow \kF\}}{\#\text{Aut}(\kL)}
=v^{-\langle (r-1,d-N),(1,N)\rangle}\sum_{\kL \in Pic^N(X)}\frac{\#\{\kL \hookrightarrow \kF\}}{v^{-2}-1}.
\end{equation*}
Let us decompose $\kF=\mathcal{V}_{\kF}\oplus \mathcal{T}_{\kF}$ into a direct sum of a vector bundle and a torsion sheaf, and let us assume that $\kF \in \textbf{C}_{\geq n}$. Then $\kF \in \textbf{C}_{\geq 2(g-1)+N}$ and thus
$\text{Ext}(\kL,\kF)=0$ by Serre duality. This in turn implies that $dim\;\text{Hom}(\kL,\kF)=\langle (1,N),(r,d)\rangle$. Any nonzero map from a line bundle to a vector bundle is an embedding. From this we deduce that
\begin{equation*}
\#\{\kL \hookrightarrow \kF\}=v^{-2dim\;\text{Hom}(\kL,\kF)}-v^{-2dim\;\text{Hom}(\kL,\mathcal{T}_{\kF})}
=v^{-2\langle (1,N),(r,d)\rangle}-v^{-2deg(\mathcal{T}_{\kF})}.
\end{equation*}
The important point is that this only depends on $deg(\mathcal{T}_{\kF})$. From this discussion we deduce that there exists nonzero constants $c_l$ for $l \geq 0$ such that
$$ \mathbf{1}_{r-1,d-N}\cdot\mathbf{1}^{\textbf{vec}}_{1,N}\in c_0 \mathbf{1}^{\textbf{vec}}_{r,d}+\sum_{l = 1}^{d-rn} c_{l}\mathbf{1}^{\textbf{vec}}_{r,d-l}\cdot\mathbf{1}_{0,l} + \H_{\E}^{< n}.$$
We may rewrite this last equation as
$$c_0\mathbf{1}^{\textbf{vec}}_{r,d} \in \mathbf{1}_{r-1,d-N}\cdot\mathbf{1}^{\textbf{vec}}_{1,N}-\sum_{l = 1}^{d-rn} c_{l}\mathbf{1}^{\textbf{vec}}_{r,d-l}\cdot\mathbf{1}_{0,l} + \H_{\E}^{< n}.$$
Now, by our two induction hypotheses we have $\mathbf{1}_{r-1,d-N} \in \widehat{\UU}_{\E}$ and
$\mathbf{1}^{\textbf{vec}}_{r,d-l} \in \widehat{\UU}_{\E}$ for all $l \geq 1$. But then (\ref{E:huyt}) holds as well.
We are done. \end{proof}

\vspace{.1in}

\noindent
\textit{Proof of Theorem~\ref{T:sstable}.} We have to show that $\mathbf{1}^{\ss}_{\a}$ belongs to $\UU_{\E}$, and not only to $\widehat{\UU}_{\E}$.
By Proposition~\ref{P:A4}, there exists for all $n$ an element $v_n \in \H^{<n}_{\E}$ such that $u_n:= \mathbf{1}^{\ss}_{\a} + v_n \in \UU_{\E}$. We may further decompose
$v_n=\sum_{\underline{\a}}v_{n,\underline{\a}}$ according to the HN type $\underline{\a}$. The set of $\underline{\a}$ for which $v_{n,\underline{\a}}$ is nonzero is finite since $v_n \in \H_{\E}$. Our proof is based on the following two lemmas.

\vspace{.1in}

\begin{lem}\label{L:A1} There exists $n \ll 0$ such that for any HN type $\underline{\a}=(\a_1, \ldots, \a_l)$ of weight $\alpha$ satisfying $\mu(a_1)< n$, we have $\mu(\a_{i+1})-\mu(\a_i) > 2(g-1)$ for some $1 \leq i \leq l$.
\end{lem}
\begin{proof} Let $\underline{\a}=(\a_1, \ldots, \a_l)$ be as above. We have $deg(\a)=rank(\a_1)\mu(\a_1) + \cdots +_ rank(\a_l) \mu(\a_l)$.
If $\mu(\a_1)<n$ and $\mu(\a_{i+1})-\mu(\a_i) \leq 2(g-1)$ for all $i$ then
\begin{equation*}
\begin{split}
deg(\a) &< rank(\a_1)n + rk(\a_2)(n+2(g-1)) + \cdots + rank(\a_l)(n+2(g-1)(l-1)\\
&=rank(\a)n +\sum_{i=2}^{l} 2(g-1)(l-1)rank(\a_i)\\
&< rank(\a) \bigg(n + 2(g-1)\sum_{l=1}^{rank(\a)} l\bigg).
\end{split}
\end{equation*}
This is impossible for $n$ sufficiently negative.
\end{proof}

\vspace{.1in}

\begin{lem}\label{L:A2} Let $\mathcal{F} \in Coh(\E)$ be a coherent sheaf of class $\alpha$ and of HN type $(\a_1, \ldots, \a_l)$.
Assume that $\mu(\a_{i+1})-\mu(\a_i) > 2(g-1)$ for some $i$. Then $1_{\mathcal{F}}=m \circ \Delta_{\beta,\gamma}(1_{\mathcal{F}})$ for
$\beta=\a_1 + \cdots + \a_i, \gamma=\a_{i+1} + \cdots + \a_l$.
\end{lem}
\begin{proof} Let $\mathcal{F}_l \subset \cdots \subset \mathcal{F}_1=\mathcal{F}$ be the HN filtration of $\mathcal{F}$. Since $\mathcal{F}_{i+1} \in \mathbf{C}_{\geq \mu(\a_{i+1})}$ and $\mathcal{F}/\mathcal{F}_{i+1} \in \mathbf{C}_{\leq \mu(a_i)}$ while $\mu(\a_{i+1})-\mu(\a_i) > 2(g-1)$ we have $\text{Ext}(\mathcal{F}_{i+1}, \mathcal{F}/\mathcal{F}_{i+1})=0$ (see Proposition~\ref{P:A1}). It follows that $\mathcal{F} \simeq \mathcal{F}_{i+1} \oplus \mathcal{F}/\mathcal{F}_{i+1}$. Moreover, $1_{\mathcal{F}/\mathcal{F}_{i+1}} 1_{\mathcal{F}_{i+1}}=v^{-\langle \mathcal{F}/\mathcal{F}_{i+1}, \mathcal{F}_{i+1}\rangle} 1_{\mathcal{F}}$ since there is a unique subsheaf of $\mathcal{F}$ isomorphic to $\mathcal{F}_i$. Hence Lemma~\ref{L:A2} will be proved once we show that $\Delta_{\beta,\gamma}(1_{\mathcal{F}})=v^{\langle \mathcal{F}/\mathcal{F}_{i+1}, \mathcal{F}_{i+1}\rangle}1_{\mathcal{F}/\mathcal{F}_{i+1}} \otimes 1_{\mathcal{F}_{i+1}}$. But this last equation is a consequence of the fact that
there exists a \textit{unique} subsheaf of $\mathcal{F}$ of class $\gamma$, namely $\mathcal{F}_{i+1}$ (see Proposition~\ref{P:A2}).
\end{proof}

\vspace{.1in}

We are now ready to finish the proof of Theorem~\ref{T:sstable}. Let us choose some $n \ll 0$ as in Lemma~\ref{L:A1}. Let
$\mathcal{A}$ be the (finite) set of all $\underline{\a}$ for which $v_{n, \underline{\a}}$ is nonzero and let $\underline{\a}^0$ be the lower boundary
of the convex hull of elements of $\mathcal{A}$.

\vspace{.2in}

\centerline{
\begin{picture}(200,120)
\put(0,0){\line(0,1){120}}
\put(0,60){\line(1,0){200}}
\put(0,60){\circle*{2}}
\put(-7,60){$\textbf{o}$}
\put(10,30){\circle*{2}}
\put(120,20){\circle*{2}}
\put(10,30){\line(4,1){120}}
\put(130,60){\line(5,4){50}}
\put(130,60){\circle*{2}}
\put(182,98){$\a$}
\put(60,10){\circle*{2}}
\put(0,60){\line(6,-5){60}}
\put(180,100){\circle*{2}}
\put(0,60){\line(3,-1){120}}
\thicklines
\put(0,60){\line(1,-3){10}}
\put(120,20){\line(3,4){60}}
\put(60,10){\line(6,1){60}}
\put(10,30){\line(5,-2){50}}
\thinlines
\end{picture}}
\vspace{.05in}
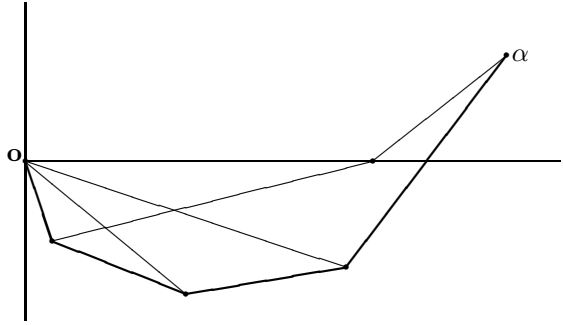
\centerline{\textbf{Figure 2.} The convex hull of a set of HN polygons.}

\vspace{.15in}

Thus $\underline{\a}^0=(\a_1^0, \ldots, \a_m^0)$ is also a convex path in $\Z^2$ of
weight $\a$. Moreover $\mu(\a_1^0)<n$ so that the conclusion of Lemma~\ref{L:A1} applies. Choose $i$ such that $\mu(\a_{i+1}^0)-\mu(\a_i^0) > 2(g-1)$ and set $\beta=\a_1^0 + \cdots + \a_i^0, \gamma=\a_{i+1}^0 + \cdots + \a_m^0$. By Lemma~\ref{L:A2}, $\Delta_{\beta,\gamma}(1_{\mathcal{F}})=0$ for all sheaves $\mathcal{F}$ whose HN polygon doesn't lie below the segment $\beta$.

\vspace{.2in}

\centerline{
\begin{picture}(200,120)
\put(0,0){\line(0,1){120}}
\put(0,60){\line(1,0){200}}
\put(0,60){\circle*{2}}
\put(-7,60){$\textbf{o}$}
\put(10,30){\circle*{2}}
\put(120,20){\circle*{2}}
\put(10,30){\line(4,1){120}}
\put(130,60){\line(5,4){50}}
\put(130,60){\circle*{2}}
\put(182,98){$\a$}
\put(60,10){\circle*{2}}
\put(57,0){$\beta$}
\put(0,60){\line(6,-5){60}}
\put(180,100){\circle*{2}}
\put(0,60){\line(3,-1){120}}
\thicklines
\put(0,60){\line(1,-3){10}}
\put(120,20){\line(3,4){60}}
\put(60,10){\line(6,1){60}}
\put(10,30){\line(5,-2){50}}
\thinlines
\end{picture}}
\vspace{.05in}
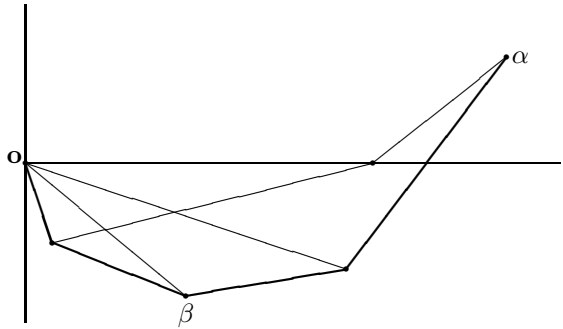
\centerline{\textbf{Figure 3.} Choice of the vertex $\beta$.}

\vspace{.15in}

This implies that $\Delta_{\beta,\gamma}(v_{n,\underline{\a}})=0$ for all HN types $\underline{\a}$ whose associated polygon doesn't pass through the point $\beta$. Furthermore, by Lemma~\ref{L:A2} again, $m \circ \Delta_{\beta,\gamma}(v_{n, \underline{\a}})=v_{n,\underline{\a}}$ for any HN type $\underline{\a}$ whose polygon does pass through $\beta$. Hence
$$m \circ \Delta_{\beta,\gamma}(u_n)=m \circ \Delta_{\beta,\gamma}\bigg(\mathbf{1}^{\ss}_{\a} + \sum_{\underline{\a}} v_{n,\underline{\a}}\bigg)=\sum_{\underline{\a} \in Z_{\beta}} v_{n,\underline{\a}},$$
where $Z_{\beta}$ is the set of all HN types passing through $\beta$. Because $u_n$ belongs to $\UU_{\E}$, which is stable under the coproduct,
we deduce that $\sum_{\underline{\a} \in Z_{\beta}} v_{n,\underline{\a}}$ belongs to $\UU_{\E}$ as well. Hence the same holds for $u'_n=\mathbf{1}^{\ss}_{\a} + \sum_{\underline{\a} \notin Z_{\beta}} v_{n,\underline{\a}}$. Notice that $u'_n$ contains strictly fewer terms than $u_n$. Arguing as above repeatedly we obtain better and better approximations of $\mathbf{1}^{\ss}_{\a}$ by elements of $\UU_{\E}$ until we arrive at $\mathbf{1}^{\ss}_{\a} \in \UU_{\E}$. Theorem~\ref{T:sstable} is proved.\qed

\vspace{.1in}

The combination of Theorem~\ref{T:sstable} and Proposition~\ref{P:A3} yields the following result~:

\vspace{.1in}

\begin{cor} For any HN type $\underline{\a}$ we have $\mathbf{1}_{S_{\underline{\a}}} \in \UU_{\E}$.
\end{cor}

\vspace{.2in}

\addtocounter{theo}{1}

\paragraph{\textbf{Remark \thetheo}} The above proof actually shows that $\widehat{\UU}_{\E} \cap \H_{\E}=\UU_{\E}$.  

\vspace{.2in}

\section{Spherical Eisenstein sheaves}

\vspace{.15in}

\paragraph{\textbf{4.1.}} Let us set $\overline{X}=X \otimes \fqb$. For $\a \in K'_0(\overline{X})=\Z^2$, let $\underline{Coh}^\a_X$
stand for the moduli stack parametrizing coherent sheaves of class $\a$ over $\overline{X}$. This is a smooth irreducible stack, which is locally of finite type (see e.g. \cite{Laumon} ). It carries a Harder-Narasimhan stratification $\underline{Coh}^\a_X=\bigsqcup_{\underline{\a}} \underline{S}_{\underline{\a}}$ similar to the one existing for $\mathcal{I}_X$, and each locally closed substack $\underline{S}_{\underline{\a}}$ is smooth and of finite type.

\vspace{.1in}

We will now define, following \cite{SInvent}, a certain category of simple perverse sheaves
over the stacks $\underline{Coh}^\a_X$. For this we consider the following induction diagrams, for $\a,\beta \in K'_0(X)$~:

\begin{equation}\label{E:Inddiagram}
\xymatrix{ & \underline{\mathcal{E}}_{\a,\beta}  \ar[dl]_-{p_1} \ar[dr]^-{p_2} &\\
\underline{Coh}^{\a}_X \times \underline{Coh}^{\beta}_X & & \underline{Coh}^{\a+\beta}_X}
\end{equation}
where $\underline{\mathcal{E}}_{\a,\beta}$ is the stack classifying inclusions $\mathcal{G} \subset \mathcal{F}$ of a coherent sheaf
$\mathcal{G}$ over $\overline{X}$ of class $\beta$ into a coherent sheaf $\mathcal{F}$ over $\overline{X}$ of class $\a +\beta$; and where the maps $p_1$, $p_2$ are given by the functors $\mathcal{G} \subset \mathcal{F} \mapsto (\mathcal{F}/\mathcal{G}, \mathcal{G})$ and $\mathcal{G} \subset \mathcal{F} \mapsto \mathcal{F}$. The morphism $p_1$ is smooth while $p_2$ is proper and representable (see \cite{Laumon}).

\vspace{.1in}

Let $D^b(\underline{Coh}^\a_X)$ be the bounded derived category of constructible $\qlb$-sheaves over $\underline{Coh}^\a_X$. We define induction and restriction functors as
\begin{equation}\label{E:killbill}
\begin{split}
\underline{m}~: D^b(\underline{Coh}^{\alpha}_X \times \underline{Coh}^{\beta}_{X}) &\to  D^b(\underline{Coh}^{\alpha+\beta}_{X} )\\
\mathbb{P} &\mapsto p_{2!} p_1^* (\mathbb{P})[dim\;p_1],
\end{split}
\end{equation}
and 
\begin{equation}
\begin{split}
\underline{\Delta}~: D^b(\underline{Coh}^{\alpha+\beta}_{X})  &\to  D^b(\underline{Coh}^{\alpha}_{X} \times \underline{Coh}^{\beta}_{X}) \\
\mathbb{P} &\mapsto p_{1!} p_2^* (\mathbb{P})[dim\;p_2].
\end{split}
\end{equation}
By the Decomposition Theorem of \cite{BBD}, $\underline{m}$ preserves the subcategory of semisimple complexes of geometric origin.   Both of the above functors are associative in the appropriate sense. We will sometimes write $\mathbb{P} \star \mathbb{Q}$ for $\underline{m}(\mathbb{P} \boxtimes \mathbb{Q})$. For $\a \in K'_0(X)$, let
$\mathbf{1}_{\a}={\qlb}_{\underline{Coh}^{\a}_X}[dim\;  \underline{Coh}^{\a}_X]$ be the constant complex over $\underline{Coh}^\a_X$. We will call a product of the form
$$L_{\a_1, \ldots, \a_r} = \mathbf{1}_{\a_1} \star \cdots \star \mathbf{1}_{\a_r}$$
a \textit{Lusztig sheaf}. It is a semisimple complex. We let $\mathcal{P}_X=\bigsqcup_{\a} \mathcal{P}^{\a}$ stand for the set of all simple perverse sheaves appearing in some Lusztig sheaf $L_{\a_1, \ldots, \a_r}$ where for all $\a_i=(r_i,d_i)$ we have $r_i \leq 1$. We denote by $\mathcal{Q}_X=\bigsqcup_{\a} \mathcal{Q}^{\alpha}$ the additive category generated by the objects of $\mathcal{P}_X$ and their shifts. 

\vspace{.2in}

\section{IC sheaves of Harder-Narasimhan strata}

\vspace{.15in}

The purpose of this section is to prove the following result~:

\vspace{.1in}

\begin{theo}\label{T:maingeom} For any Harder-Narasimhan type $\underline{\a}$ we have $IC(\underline{S}_{\underline{\a}}) \in \mathcal{P}_X$.
\end{theo}

\vspace{.1in}

This can be viewed as a direct geometric analog of Theorem~\ref{T:sstable}. We will first establish the following special case~:

\vspace{.1in}

\begin{prop}\label{P:stablegeom} For any $\a \in K'_0(X)$ we have $\mathbf{1}_{\a} \in \mathcal{P}_X$.
\end{prop}
\begin{proof} We argue by induction on the rank $r$ of $\a$. If $r=1$ then $\a=(1,d)$ for some $d$ and by definition we have $\mathbf{1}_{\underline{Coh}^{(1,d)}_X} \in \mathcal{P}_X$. Let us fix some $\a$ of rank $r >1$ and  let us assume that $\mathbf{1}_{\beta}$ belongs to $\mathcal{P}_X$ for all $\beta$ of rank strictly less than $r$. Let us choose some $d \ll \mu(\a)$ and set $\beta=\a-(1,d)$. Consider the convolution diagram
\begin{equation}\label{E:Inddiagram2}
\xymatrix{ & \underline{\mathcal{E}}_{\beta,(1,d)}  \ar[dl]_-{p_1} \ar[dr]^-{p_2} &\\
\underline{Coh}^{\beta}_X \times \underline{Coh}^{1,d}_X & & \underline{Coh}_X^{\a}}
\end{equation}
and the corresponding product $\mathbf{1}_{\beta} \star \mathbf{1}_{(1,d)} = p_{2!} (\mathbf{1}_{\underline{\mathcal{E}}_{\beta,(1,d)}})$. By the induction hypothesis, both $\mathbf{1}_{\beta}$ and $\mathbf{1}_{(1,d)}$ belong to $\mathcal{P}_X$. Hence it is enough to show that $Hom(\mathbf{1}_{\a}, \mathbf{1}_{\beta} \star \mathbf{1}_{(1,d)}) \neq \{0\}$. The strata $\underline{S}_{\a} \subset \underline{Coh}^\a_X$ is open and we have $\mathbf{1}_{\a}=IC(\underline{S}_{\a})$. Because the complex $\mathbf{1}_{\beta} \star \mathbf{1}_{(1,d)}$ is semisimple, it is sufficient to prove that 
\begin{equation}\label{E:poliu}
Hom\big(\mathbf{1}_{\underline{S}_{\a}}, (\mathbf{1}_{\beta} \star \mathbf{1}_{(1,d)})_{|\underline{S}_{\a}}\big)=Hom\big(\mathbf{1}_{\underline{S}_\a}, p_{2!}(\mathbf{1}_{\underline{\mathcal{E}}_{\beta,(1,d)}})_{|\underline{S}_{\a}}\big) \neq \{0\}.
\end{equation}
Consider the cartesian diagram obtained by restricting (\ref{E:Inddiagram2}) to the open strata $\underline{S}_{\a}$
\begin{equation}
\xymatrix{ p_2^{-1}(\underline{S}_{\a}) \ar[r]^{j} \ar[d]_-{p'_2} & \underline{\mathcal{E}}_{\beta,(1,d)} \ar[d]^{p_2} \\
\underline{S}_{\a}  \ar[r]^{j'} & \underline{Coh}^\a_X}
\end{equation}
Here $j, j'$ are the open embeddings. We now describe the map $p'_2: p_2^{-1}(\underline{S}_{\a}) \to \underline{S}_{\a}$ explicitly.
Observe that any semistable sheaf $\mathcal{F}$ of class $\a$ is a vector bundle, and hence any subsheaf $\mathcal{G} \subset \mathcal{F}$ of class $(1,d)$ is a line bundle. Let $\underline{\mathcal{H}}$ be the $Hom$-stack over $\underline{Pic}^d \overline{X} \times \underline{S}_{\a}$, i.e. the stack parametrizing triples
$ (\mathcal{F},\mathcal{G},f)$ where $\mathcal{F} \in S_{\a} \overline{X}, \;\mathcal{G} \in Pic^d \overline{X}$ and $f \in Hom(\mathcal{G},\mathcal{F})$. Let $a: \underline{\mathcal{H}} \to \underline{Pic}^d \overline{X} \times \underline{S}_{\a}$ be the projection.
Observe that since $d \ll \a$ and since $\mathcal{F}$ is semistable, we have $dim\; Hom(\mathcal{G}, \mathcal{F}) = dim\; Hom(\mathcal{G},\mathcal{F})-dim\;Ext(\mathcal{G},\mathcal{F}) =\langle (1,d), \a\rangle$ for any $\mathcal{F},\mathcal{G}$ as above. It follows that
$a$ is a vector bundle of rank $\langle (1,d), \a\rangle$. Denote by $a'~: \underline{\mathcal{H}}' \to \underline{Pic}^d \overline{X} \times \underline{S}_{\a}$ the associated projective bundle. It is easy to see that $p'_2 :p_2^{-1}(\underline{S}_{\a}) \to \underline{S}_{\a}$ is canonically isomorphic to the composition of $a'$ with the projection $\pi: \underline{Pic}^d \overline{X} \times \underline{S}_{\a} \to \underline{S}_{\a}$.

We may now compute
$$(\mathbf{1}_{\beta} \star \mathbf{1}_{(1,d)})_{|\underline{S}_{\a}}=(j')^* p_{2!} (\mathbf{1}_{\underline{\mathcal{E}}_{\beta,(1,d)}})
= p'_{2!}( \mathbf{1}_{p_2^{-1}(\underline{S}_{\a})})=\pi_! a'_! (\mathbf{1}_{\underline{\mathcal{H}}'}).$$
Because $a'$ is a projective bundle, the constant sheaf $\mathbf{1}_{\underline{Pic}^d \overline{X} \times \underline{S}_{\a}}$ appears as a direct summand of $a'_! (\mathbf{1}_{\underline{\mathcal{H}}'})$. Therefore $\mathbf{1}_{\underline{S}_{\a}}$ appears as a direct
summand of $\pi_! a'_! (\mathbf{1}_{\underline{\mathcal{H}}'}).$ This shows (\ref{E:poliu}) and finishes the proof of Proposition~\ref{P:stablegeom}.
\end{proof}

\vspace{.15in}

We are now in position to prove Theorem~\ref{T:maingeom}. Let $\underline{\a}=(\a_1, \ldots, \a_l)$ be some HN type of weight $\a=\sum \a_i$. By Proposition~\ref{P:stablegeom}, all the perverse sheaves $IC(\underline{S}_{\a_i})=\mathbf{1}_{\a_i}$ belong to $\mathcal{P}_X$. We will show that $IC(\underline{S}_{\underline{\a}})$ belongs to $\mathcal{P}_X$ as well by proving that
\begin{equation}\label{E:poliu2}
Hom\big(IC(\underline{S}_{\underline{\a}}), \mathbf{1}_{\a_1} \star \cdots \star \mathbf{1}_{\a_l}\big) \neq \{0\}.
\end{equation}
For this, consider the (iterated) induction diagram
\begin{equation}\label{E:Inddiagram3}
\xymatrix{ & \underline{\mathcal{E}}^{(l)}_{\a_1,\ldots, \a_l}  \ar[dl]_-{p_1} \ar[dr]^-{p_2} &\\
\underline{Coh}^{\a_1}_X \times \cdots \times \underline{Coh}^{\a_l}_X & & \underline{Coh}_X^{\a}}
\end{equation}

We claim that $\underline{S}_{\underline{\a}}$ is open and dense in $Im(p_2)$. Indeed, set $\underline{\mathcal{E}}^{(l),0}_{\a_1,\ldots, \a_l}=p_{1}^{-1}(\underline{S}_{\a_1} \times \cdots \underline{S}_{\a_l})$. It is an open dense subset of $\underline{\mathcal{E}}^{(l)}_{\a_1,\ldots, \a_l}$ and by construction $p_2(\underline{\mathcal{E}}^{(l),0}_{\a_1,\ldots, \a_l})=\underline{S}_{\underline{\a}}$. Since $p_2$ is continuous, $p_2^{-1}( \underline{Coh}^\a_X \backslash \overline{\underline{S}_{\underline{\a}}})$ is an open substack of
$\underline{\mathcal{E}}^{(l)}_{\a_1,\ldots, \a_l}$ which does not intersect $\underline{\mathcal{E}}^{(l),0}_{\a_1,\ldots, \a_l}$. But this means that $p_2^{-1}( \underline{Coh}^\a_X \backslash \overline{\underline{S}_{\underline{\a}}})$ is empty, i.e. that $Im(p_2) \subset 
 \overline{\underline{S}_{\underline{\a}}}$ as wanted.
 
\vspace{.1in}
 
By definition, $\mathbf{1}_{\a_1} \star \cdots \star \mathbf{1}_{\a_l} = p_{2!} (\mathbf{1}_{\underline{\mathcal{E}}^{(l)}_{\a_1, \ldots, \a_l}})$. This is a semisimple complex and, by the above, $\underline{S}_{\underline{\a}}$ is open in its support. Therefore 
\begin{equation*}
Hom\big(IC(\underline{S}_{\underline{\a}}), \mathbf{1}_{\a_1} \star \cdots \star \mathbf{1}_{\a_l}\big)=Hom\big(\mathbf{1}_{\underline{S}_{\underline{\a}}}, j_{\underline{\a}}^*( \mathbf{1}_{\a_1} \star \cdots \mathbf{1}_{\a_l})\big)=Hom\big(\mathbf{1}_{\underline{S}_{\underline{\a}}}, j_{\underline{\a}}^*p_{2!}(\mathbf{1}_{\underline{\mathcal{E}}^{(l)}_{\a_1, \ldots, \a_l}})\big)
\end{equation*}
where $j_{\underline{\a}} : \underline{S}_{\underline{\a}} \to \underline{Coh}^{\a}_X$ denote the inclusion. Observe that by the uniqueness of the Harder-Narasimhan filtration $\mathcal{F}_{l} \subset \cdots \subset \mathcal{F}_1= \mathcal{F}$ of a coherent sheaf $\mathcal{F} \in S_{\underline{\a}}$, the projective map $p_2$ restricts to an isomorphism $p_2^{-1}(\underline{S}_{\underline{\a}}) \stackrel{\sim}{\to} \underline{S}_{\underline{\a}}$. By base change, $j_{\underline{\a}}^*p_{2!}(\mathbf{1}_{\underline{\mathcal{E}}^{(l)}_{\a_1, \ldots, \a_l}})=\mathbf{1}_{\underline{S}_{\underline{\a}}}$. But then $Hom\big(IC(\underline{S}_{\underline{\a}}), \mathbf{1}_{\a_1} \star \cdots \star \mathbf{1}_{\a_l}\big) =\qlb$, and (\ref{E:poliu2}) follows. Theorem~\ref{T:maingeom} is proved.
\qed

\vspace{.2in}

\addtocounter{theo}{1}

\paragraph{\textbf{Remark \thetheo}} The collection of simple perverse sheaves $\{IC(\underline{S}_{\underline{\a}}\}_{\underline{\a}}$ by no means exhausts of all $\mathcal{Q}_X$. We illustrate this by a very simple example showing that one has \textit{at least} to consider
the various Brill-Noether stratas $\underline{W}^k_{r,d}$ of the stacks $\underline{Coh}^{r,d}_X$ (see e.g. \cite{ACGH} or \cite{King}; note that we use slightly different notation). Let $(r,d)=(2,0)$, and let us consider the direct summands of
$\mathbf{1}_{(1,0)} \star \mathbf{1}_{(1,0)}$. The stack $\underline{S}_{2,0} \subset \underline{Coh}^{2,0}_X$ of semistable bundles may be stratified as follows~: $\underline{S}_{2,0}=\underline{W}_{2,0}^{0} \cup \underline{W}^1_{(2,0)} \cup \underline{W}^2_{2,0} \cup \underline{U}$ where
$$\underline{W}_{2,0}^0=\{ \mathcal{V} \in {S}_{2,0}\;|\; Hom(\mathcal{L}, \mathcal{V})=\{0\}, \;\forall\; \mathcal{L} \in Pic^0 \overline{X}\},$$
$$\underline{W}_{2,0}^1=\{ \mathcal{V} \in {S}_{2,0}\;|\; \exists\; ! \;\mathcal{L} \in Pic^0 \overline{X},\;Hom(\mathcal{L}, \mathcal{V})=\fqb\},$$
$$\underline{U}=\{  \mathcal{V} \in {S}_{2,0}\;|\; \exists\; \mathcal{L}, \mathcal{L}' \in Pic^0\overline{X},\; \mathcal{L}\neq \mathcal{L}',\; \mathcal{V} \simeq \mathcal{L} \oplus \mathcal{L}'\},$$
$$\underline{W}_{2,0}^2=\{ \mathcal{V} \in {S}_{2,0}\;|\; \exists\; ! \;\mathcal{L} \in Pic^0 \overline{X},\;Hom(\mathcal{L}, \mathcal{V})=\fqb^2\}.$$
The strata $\underline{W}^1_{2,0}$ consists of semistable vector bundles which are nontrivial extensions 
$$\xymatrix{ 0 \ar[r] & \mathcal{L} \ar[r] & \mathcal{V} \ar[r] & \mathcal{L}' \ar[r] & 0}$$ 
of two degree zero line bundles $\mathcal{L}, \mathcal{L}'$; the strata $\underline{W}_{2,0}^2$ consists of semistable bundles of the form $\mathcal{V} \simeq \mathcal{L} ^{\oplus 2}$ for some degre zero line bundle $\mathcal{L}$. Moreover, $\underline{W}^0_{2,0}$ is open dense, $\underline{W}^2_{2,0}$ is closed and we have inclusions of strata closures $\overline{\underline{W}^1_{2,0}} \supset \overline{\underline{U}} \supset {\underline{W}^2_{2,0}}$.

\vspace{.1in}

The restriction 
$$p'_2~: p_2^{-1}(\underline{S}_{2,0}) \to \underline{S}_{2,0}$$
of the proper map $p_2$ in the induction diagram (\ref{E:Inddiagram}) corresponding to $\mathbf{1}_{(1,0)} \star \mathbf{1}_{(1,0)}$ respects the above stratification. The following table lists the dimension of each strata as well as the type of fiber~:

$$
\begin{tabular}{|c|c|c|}
\hline
Strata & Dimension & Fiber \\
\hline
$\underline{W}^0_{2,0}$ & $4g-4$ & $\emptyset$\\
\hline
$\underline{W}^1_{2,0}$ & $3g-3$ & $\{pt\}$\\
\hline
$\underline{U}$& $2g-2$ & $\{pt\} \cup \{pt\}$\\
\hline
$\underline{W}^2_{2,0}$& $g-4$ & $\mathbb{P}^1$\\
\hline
\end{tabular}
$$

\vspace{.1in}

It follows from the above table that $p'_2$ is a small resolution of the closure $\overline{\underline{W}^1_{2,0}}$ of $\underline{W}^1_{2,0}$, and hence that $p'_{2!}(\mathbf{1}_{p_2^{-1}(\underline{S}_{2,0})}) =IC(\underline{W}^1_{2,0})$. Considering induction products of the form $\mathbf{1}_{1,d} \star \mathbf{1}_{1,-d}$ for $d=1, \ldots, 2g-2$ we obtain elements of $\mathcal{P}_X$ supported on other nontrivial Brill-Noether type stratas of $\underline{S}_{2,0}$.

\vspace{.2in}

\addtocounter{theo}{1}

\paragraph{\textbf{Remark \thetheo}} The method of proof of Theorem~\ref{T:maingeom} is readily transposable to the context of
quivers. Let $\vec{Q}$ be a quiver and let us assume that it contains no oriented cycles. Lusztig defined in \cite{Lusztig} a set $\mathcal{Q}_{\vec{Q}}$ of simple perverse sheaves over the moduli stacks $\underline{\mathcal{M}}^{\a}_{\vec{Q}}$, $\a \in K_0(\vec{Q})$.
To any linear form $\Theta : K_0(\vec{Q}) \to \mathbb{C}$ (the `stability parameter') one may attach a slope function $\mu_{\Theta}$ on $K_0(\vec{Q})$ and a Harder-Narasimhan stratification $\underline{\mathcal{M}}^{\a}_{\vec{Q}}=\bigsqcup_{\underline{\a}} \underline{S}_{\underline{\a}}$ (see \cite{Reineke}). Then

\vspace{.05in}

\begin{theo} For any stability parameter $\Theta$ and any Harder-Narasimhan strata $\underline{S}_{\underline{\a}}$ we have $IC(\underline{S}_{\underline{\a}}) \in \mathcal{P}_{\vec{Q}}$.
\end{theo}

Note that the analog of Proposition~\ref{P:stablegeom} holds since $\vec{Q}$ has no oriented cycles. 
Theorem~\ref{T:sstable} also holds in the context of quivers, where it is a simple consequence of Reineke's inversion formula (\cite{Reineke}).

\vspace{.2in}

\small{}

\vspace{4mm}

\noindent
O. Schiffmann, \texttt{olive@math.jussieu.fr},\\
D\'epartement de Math\'ematiques, Universit\'e de Paris 6, 175 rue du Chevaleret, 75013 Paris, FRANCE.

\vspace{.1in}

\end{document}